\documentclass[12pt]{article}
\usepackage[mathscr]{eucal}
\usepackage{amssymb,latexsym}
\usepackage{verbatim}
\usepackage{amsmath}
\usepackage{amsthm}
\usepackage{enumerate}
\usepackage{authblk}
\usepackage{color}

\newtheorem{thm}{Theorem}[section]
\newtheorem{lem}[thm]{Lemma}
\newtheorem{Def}[thm]{Definition}
\newtheorem{prop}[thm]{Proposition}
\newtheorem{cor}[thm]{Corollary}
\newtheorem{ex}[thm]{Example}

\newcommand\vs{\vspace}

\newcommand\ol{\overline}

\newcommand\wt{\widetilde}

\renewcommand\S{\Sigma}

\newcommand\s{\sigma}
\renewcommand\d{\partial}

\newcommand\g{\gamma}

\renewcommand\t{\tau}

\newcommand\beq{\begin{equation}}
\newcommand\eeq{\end{equation}}
\newcommand\ben{\begin{enumerate}}
\newcommand\een{\end{enumerate}}
\newcommand\bit{\begin{itemize}}
\newcommand\eit{\end{itemize}}

\newcommand{\R}{\mathbb R}

\newcommand{\ov}{\overline}

\newcommand{\ext}{\text{{\rm ext}}}
\newcommand{\pd}{\partial}

\newcounter{mnotecount}

\title{Some Remarks on the $C^0$-(in)extendibility of Spacetimes}

\author{Gregory J. Galloway
}
\author{Eric Ling}
\affil{Department of Mathematics
\\ University of Miami }

\begin{document}
\date{}
\maketitle
\vspace{.2in}

\begin{abstract} 

The existence, established over the past number of years and supporting earlier work of Ori \cite{Ori},
of physically relevant black hole spacetimes that admit $C^0$ metric extensions beyond the future Cauchy horizon, while being $C^2$-inextendible,  has focused attention on fundamental issues concerning the strong cosmic censorship conjecture.  These issues were recently discussed in the work of Jan Sbierski \cite{Sbierski}, in which he established the (nonobvious) fact that the Schwarschild solution in global Kruskal-Szekeres coordinates is $C^0$-inextendible.  In this paper we review aspects of Sbierski's methodology in a general context, and use similar techniques, along with some new observations, to consider the $C^0$-inextendibility 
of open FLRW cosmological models. We find that a certain special class of open FLRW spacetimes, which we have dubbed `Milne-like,' actually admit $C^0$ extensions through the big bang. For spacetimes that are not Milne-like, we prove some inextendibility results within the class of spherically symmetric spacetimes.

\end{abstract}


\section{Introduction}

Ever since the realization, from the Hawking-Penrose singularity theorems, that singularities in spacetime can develop under generic circumstances, the question has been considered as to what extent General Relativity is a classically deterministic theory.  The essence of Penrose's strong cosmic censorship conjecture is that, indeed, GR is deterministic.  Put in rough physical terms, under reasonable physical conditions, spacetime should not develop {\it naked singularities}, that is to say, no singularity (due e.g.\ to curvature blow-up) should ever be visible to any observer.  Such singularities would undermine the predictive ability of GR.  Penrose formulated the absence of naked singularities in terms of terminal indecomposable past sets (TIPs), which are defined completely in terms of the causal structure of spacetime:  Strong cosmic censorship requires that no singular TIP be contained in the timelike past of some point in spacetime; see e.g. \cite{PenroseCC}.  

More modern statements of the strong cosmic censorship conjecture focus on the Cauchy problem for the Einstein equations, along the lines of the following.

\smallskip
\noindent
{\bf Strong cosmic censorship conjecture:} {\it The maximal globally hyperbolic development of generic compact or asymptotically flat initial data for the Einstein equations (vacuum or with reasonable matter fields) is inextendible as a suitably regular Lorentzian manifold.}

\smallskip

Given the extreme complexity of the cosmic censorship problem, efforts have been made to investigate the conjecture for restricted classes of spacetimes, e.g. classes which admit certain symmetries, or which are perturbations of exact solutions; see, for example, \cite{Ring} and \cite{Dafermos} (and references therein) for important results in the cosmological setting, and the asymptotically flat setting, respectively.  

The question arises as to what one should take as `suitably regular'.  Recent advances in our understanding of the Cauchy problem for the Einstein equations at low regularity suggests that one should perhaps allow lower than $C^2$ regularity in the statement of the conjecture.  Christodoulou \cite{Christo} established a stronger form of inextendibility, namely, that of $C^0$-inextendibility, for a generic class of spherically symmetric spacetimes satisfying the Einstein-scalar field equations.  However, subsequent work of Dafermos \cite{Dafermos03,Dafermos05} demonstrated the $C^0$-extendibility of the maximal globally hyperbolic development of solutions to the spherically symmetric Einstein-Maxwell-scalar field system, arising from small perturbations of  Reissner-Nordstrom initial data.  Moreover, more recently,  Dafermos and Luk have announced a proof, without symmetry assumptions, of the $C^0$ stability of the Kerr Cauchy horizon.   The current suggestion for the statement of the strong cosmic censorship conjecture is to require   inextendibility as a Lorentzian manifold with Christoffel symbols locally in~$L^2$.

Prior to recent work of Sbierski \cite{Sbierski}, very little had been done to address the issue of the extendibility (or not) of Lorentzian manifolds at lower regularity.   In \cite{Sbierski} Sbierski develops methods 
for establishing the $C^0$ (metric) inextendibility of Lorentzian manifolds, which he uses to prove the 
$C^0$-inextendibility of Minkowski space and the extended Schwarzchild solution.   Among several problems he lists at the end of the introduction to his paper, he poses the problem of investigating the 
$C^0$-inextendibility of cosmological spacetimes with big bang singularities, such as the FLRW models, for which some of his methods are not directly applicable.  In this paper we address some aspects of this problem.   In Section \ref{Sb} we review, in a general setting, aspects of Sbierski's methodology.  As an illustration, we demonstrate the $C^0$-inextendibility of anti-de Sitter space.  We also obtain a result concerning the structure  of the past boundary $\pd^-M$ of a spacetime $(M,g)$ within a $C^0$-extension.   In 
Section~\ref{FLRW} we introduce a general notion of `open FLRW models', and present several types of $C^0$-extensions.  In the other direction, in Section \ref{sphsym}, we present results establishing the nonexistence of spherically symmetric  $C^0$ (or $C^1$) extensions of such models.  In the last section we discuss, in a general framework, some more detailed structural properties of the past boundary $\pd^-M$ of spacetime within a $C^0$-extension, and some restrictions this imposes on the extension.

\medskip
\noindent
{\it Acknowledgments.}  The authors would like to express their thanks to Piotr Chru\'sciel for his interest in this work and for many valuable comments.   The work of GJG was partially supported  by NSF grant DMS-1313724.

\section{Sbierski's Methodology and Related Results}\label{Sb}

In this section we review, in a general spacetime setting,  Sbierki's  \cite{Sbierski} technique for  proving, for example, that Minkowski space is $C^0$-inextendible.  As an illustration, we will establish the $C^0$-inextendibilty of anti-de Sitter space.  In this framework we will also obtain a result concerning the causal structure of the boundary of spacetime within a $C^0$ extension, which we  then apply to open FLRW-type spacetimes.   This structure will be exhibited in various examples of extensions of open FLRW-type spacetimes discussed in Section \ref{FLRW}.   

\subsection{Boundaries in $C^0$-Extensions}

Manifolds will assumed to be $C^{\infty}$ and metrics will be assumed to be at least $C^0$. Let $(M^{d+1},g)$ be a spacetime, i.e.\ a connected time-oriented Lorentzian manifold. 
A spacetime $(M_\ext^{d+1}, g_\ext)$ with a continuous metric $g_\ext$ is a {\bf\emph{$C^0$-extension}} of $(M,g)$ if $(M,g)$ embeds isometrically as a proper subset of $(M_{\ext}, g_{\ext})$. We can assume that the time orientation of $(M,g)$ agrees with $(M_\ext, g_\ext)$.   
(There is no loss in generality in assuming, as we shall henceforth do, that the extended Lorentzian manifold is time-orientable.  If a non time-orientable extension did exist, then one could find a time-orientable extension by, for example, restricting the extension to a time-orientable neighborhood of a boundary point.) If no $C^0$-extension of $(M,g)$ exists, then we say that $(M,g)$ is {\bf\emph{$C^0$-inextendible}}. We denote the embedding map by $\iota: M \to M_\ext$ and henceforth, when convenient, 
identify points in $M$ with points in $\iota(M)$. Timelike curves will be piecewise smooth with right and left handed sided derivatives pointing within the same connected component of the lightcone. 

\begin{Def}
{\rm
Given a $C^0$-extension $(M_\ext,g_\ext)$ of $(M,g)$, we define:
\begin{itemize}
\item The {\bf \emph{future boundary of $M$}}, denoted by $\pd^+M$, as the set of points $p \in \pd M$ such that there exists a future directed timelike curve $\gamma: [0,1] \to M_\ext$ with $\gamma(1) = p$, $\gamma\big([0,1)\big) \subset M$.

\item The {\bf\emph{past boundary of $M$}}, denoted by $\pd^-M$, as the set of points $p \in \pd M$ such that there exists a future directed timelike curve $\gamma: [0,1] \to M_\ext$ with $\gamma(0) = p$, $\gamma\big((0,1]\big) \subset M$.
 \end{itemize}
}
\end{Def}

The following lemma gives us a way to generate points on $\pd^+M$ and $\pd^-M$. It will be used multiple times in this section.

\begin{lem}\label{inter lem}
Let $\iota: (M,g) \to(M_\ext, g_\ext)$ be a $C^0$-extension and let $\gamma: [0,1] \to M_\ext$ be a future directed timelike curve joining $p$ to $q$. 
\begin{enumerate}
\item If $p \in M$ and $q \notin M$, then $\gamma$ intersects $\pd^+M$.

\item If $p \notin M$ and $q \in M$, then $\gamma$ intersects $\pd^-M$.
\end{enumerate} 
\end{lem}

\proof
Suppose $p \in M$ and $q \notin M$. Define $s_0 = \sup \{s \in [0,1] \mid \gamma\big([0,s)\big) \subset M\}$. Since $M$ is open in $M_\ext$, we have $\gamma(s_0) \notin M$ and since $\gamma$ is future-directed, it follows that $\gamma(s_0) \in \pd^+M$. The second case is similar.
\qed

\medskip

Let $p \in M_\ext$ and $U \subset M_\ext$ be open. Recall that $I^\pm(p,U)$ are the set of points in $M_\ext$ which can be reached by a future (past) directed timelike curve whose image lies completely in $U$. If $g_\ext$ is at least $C^2$, then one can use normal neighborhoods to show that $I^\pm(p, U)$ are open in $M_\ext$. If $g$ is merely $C^0$, then $I^\pm(p,U)$ are still open, but one needs to use a different proof \cite{Sbierski, ChruGrant}. This will be needed in the following proposition which is also proved in \cite{Sbierski}.

\begin{prop}\label{bdd nonempty}
Let $\iota: (M,g) \to (M_\ext,g_\ext)$ be a $C^0$-extension. Then $\pd^+M \cup \pd^-M \neq \emptyset$. 
\end{prop}

\proof Since $M$ is a proper subset of $M_\ext$ and $M_\ext$  is connected,
$\pd M = M_\ext \setminus M$ is nonempty. Fix $p \in \pd M$ and let $U$ be an 
open neighborhood of $p$.
Fix $q \in I^-(p,U)$. By definition there exists a future directed timelike curve $\gamma: [0,1] \to U$ connecting $q$ to $p$. We either have $q \in M$ or $q \notin M$. If $q \in M$, then (1) of Lemma \ref{inter lem} implies $\pd^+M \neq \emptyset$. If $q \notin M$, then there exists an $r \in M \cap I^+(q,U)$ since $I^+(q,U)$ is an open set containing $p$. Therefore there exists a future directed timelike curve connecting $q \notin M$ to $r \in M$, so (2) of Lemma~\ref{inter lem} implies $\pd^-M \neq \emptyset$.
\qed

\medskip

Sbierski found certain sufficient conditions on $(M,g)$ to force $\pd^+M = \emptyset$. He used these to show that if Minkowski space had a $C^0$-extension then one would find $\pd^+M = \emptyset$. By reversing the time orientation, one concludes that $\pd^-M = \emptyset$. Thus the existence of a $C^0$-extension contradicts Proposition \ref{bdd nonempty}. 
Hence Minkowski space is $C^0$-inextendible. The same  conditions are used by Spierski to show that Schwarzschild spacetime is not $C^0$-extendible `beyond scri'.\footnote{In proving the 
$C^0$-inextendibility of Schwarzschild spacetime beyond the $r = 0$ singularity, Sbierski also introduced the notion of  `spacelike diameter', 
 which will not enter into the present work.}
 These  conditions are summarized in the following definition.

\medskip

\begin{Def}
{\rm
Let $(M,g)$ be a spacetime.

\vs{-.05in}
\noindent
\ben
\item[(1)] $(M,g)$ is {\bf\emph{ future one-connected}} if any two future directed timelike curves between any two points $p,q \in M$ are timelike homotopic with fixed endpoints. 

\vs{-.05in}
\item[(2)] $(M,g)$ is {\bf\emph{ future divergent}} if given any future directed inextendible timelike curve $\gamma: [0, 1) \rightarrow M$, one has $\lim_{s \rightarrow 1}d\big(\gamma(0), \gamma(s)\big) = \infty$ where $d$ is the Lorentzian distance function.
\een
} 
\end{Def}

\begin{thm}[\cite{Sbierski}]\label{future bdd empty}
Let $\iota: (M,g) \to (M_\ext, g_\ext)$ be a $C^0$-extension. If $(M,g)$ is future one-connected and future divergent then $\pd^+M = \emptyset$. 
\end{thm}

\proof[Sketch of the proof.]
Suppose there exists a point $p \in \pd^+M$. Then there exists a future directed timelike curve $\gamma:[0,1] \to M_\ext$ with $\gamma\big([0,1)\big) \subset M$ and $\gamma(1) = p$. Using the continuity of $g_\ext$ one can construct a neighborhood $U$ and a point $s_0 \in [0,1)$ such that the diamond $D = I^+\big(\gamma(s_0), U\big) \cap I^-(p,U)$ has compact closure in $U$ and finite timelike diameter. By future divergence we can find future directed timelike curves joining $\gamma(s_0)$ to $\gamma(s)$ ($s < 1$) with arbitrarily large lengths. Then future one-connectedness along with the fact that the diamond $D$ has compact closure in $U$ implies that these timelike curves must lie in the diamond $D$, so $D$ actually has infinite timelike diameter. 
\qed

\medskip

Theorem \ref{future bdd empty} enables us to obtain an  important structural result for $\pd^-M$. Recall that a set $S \subset M_\ext$ is {\bf \emph{achronal}} if for all $p,q \in S$ there exists no future directed timelike curve $\gamma: [0,1] \to M_\ext$ joining $p$ to $q$.

\begin{thm}\label{reg result}
Let $\iota: (M,g) \to (M_\ext, g_\ext)$ be a $C^0$-extension. If $\pd^+M = \emptyset$, then $\pd^-M$ is an achronal topological hypersurface.
\end{thm}

\proof We know $\pd^-M \neq \emptyset$ by Proposition \ref{bdd nonempty}. Suppose $\pd^-M$ is not achronal. Then there exist points $p,q \in \pd^-M$ with $q \in I^+(p,M_\ext)$. $I^-(q,M_\ext)$ is an open set containing $p \in \pd M$. Therefore there is a point $r \in M \cap I^-(q, M_\ext)$ and so there is a future directed timelike curve connecting $r$ to $q$. Thus $\pd^+M \neq \emptyset$ by Lemma \ref{inter lem}.

To show that $\pd^-M$ is a topological hypersurface, it suffices to show $\pd^-M \cap \text{edge}\,(\pd^-M) = \emptyset$ (\cite{HE,ON}). Fix $p \in \pd^-M$. We show $p \notin \text{edge}(\pd^-M)$. Let $r \in I^+(p,M_\ext)$ and $q \in I^-(p,M_\ext)$ and suppose $\gamma: [0,1] \to M_\ext$ is any future directed timelike curve joining $q$ to 
$r$. We need to show $\gamma$ intersects $\pd^-M$. By Lemma \ref{inter lem} it suffices to show $q \notin M$ and $r \in M$. If $q \in M$, then since $q \in I^-(p) \in M_\ext$, Lemma \ref{inter lem} implies $\pd^+M \neq \emptyset$. Since $I^-(r, M_\ext)$ is an open neighborhood of
$p$, there is a point $s \in M \cap I^-(r, M_\ext)$, and so there is a future directed timelike curve from $s$ to $r$. If $r \notin M$, then Lemma~\ref{inter lem} implies that $\pd^+M \neq \emptyset$.\qed

\medskip

\noindent{\bf Remark.}
In fact $\pd^-M$ locally satisfies a Lipschitz condition, so it has regularity~$C^{0,1}$.

\medskip

Theorems \ref{future bdd empty} and \ref{reg result} motivate   the question: which spacetimes are future one-connected and future divergent? The following subsection gives sufficient (but certainly not necessary) answers. 

\medskip

\subsection{Sufficient Conditions for Future One-Connectedness and Future Divergence}

In this subsection the spacetimes $(M,g)$ we will be interested in are warped products. Specifically $M$ will have manifold structure $M = I \times \Sigma$ where $I \subset \R$ is an open interval and $\Sigma$ is a $d$-dimensional manifold. Let $(I, dt^2)$ and $(\Sigma, h)$ be Riemannian manifolds. If $\eta: M \to I$ and $\pi: M \to \Sigma$ denote the projection maps, then the metric $g$ on $M$ is given by $g = -\eta^*dt^2 + a^2\pi^*h$ where $a: M \to (0,\infty)$ is some  smooth
function which depends only on $t$. We abuse notation and write $g = -dt^2 + a^2h$ and $a(t,p) = a(t)$.

We will show under suitable hypotheses that these spacetimes are future one-connected and future divergent. We will apply these results to open FLRW spacetimes in Section 3.

Whether or not $(M,g)$ is future one-connected depends only on its conformal class. By making the coordinate change $\tau = \int_c^t \frac{1}{a(s)}ds$, with $c \in I$, the metric becomes 
\[
g_{(\tau,p)} = a^2\big(t(\tau)\big)[-d\tau^2 + h_p].
\]

\begin{prop}\label{future connected prop}
Let $\gamma_i: [\tau_0,\tau_f] \to (M, -d\tau^2 + h)$, $i = 1,2$, be two future directed timelike curves with coinciding endpoints and each parameterized by $\tau$: $\gamma_i(\tau) = \big(\tau, \ov{\gamma}_i(\tau)\big)$. If the images of $\ov{\gamma}_1$ and $\ov{\gamma}_2$ lie completely in a common normal neighborhood $U$ of $(\Sigma,h)$ based at $\ol\g_1(\t_0) = \ol\g_2(\t_0)$, then $\g_1$ and $\g_2$ are timelike homotopic.
\end{prop}

\proof The idea of the proof is as follows (cf.\  also \cite{Sbierski}): Using the exponential map, we construct a homotopy from 
$\ov{\gamma}_1$ to the unique length minimizing geodesic connecting $\ov{\gamma}_1(\tau_0)$ and $\ov{\gamma}_1(\tau_f)$. If this homotopy is given by $\ov{\Gamma}_1(s,\tau)$, then we lift this homotopy to $M$ via $\Gamma_1(s,\tau) = 
\big(\tau, \ov{\Gamma}_1(s,\tau)\big)$, and show that $\Gamma_1$ is a timelike homotopy. We then repeat the same process for $\gamma_2$ and construct an analogous timelike homotopy $\Gamma_2$. The desired timelike homotopy is then the concatenation of $\Gamma_1$ and $\Gamma_2$. Since the procedure is symmetric, we only construct $\Gamma_1$ for $\gamma_1$ and omit the subscript. 

Let $\gamma:[\tau_0, \tau_f] \to (M, -d\tau^2 + h)$ be a future directed timelike curve with $\ov{\gamma}$ lying in a normal neighborhood  $U$ of  $(\Sigma, h)$ based at $\ol\g(\t_0)$. For each $s \in [\tau_0,\tau_f]$, let $\ov{\sigma}_s:[\tau_0, s] \to \Sigma$ be the unique length minimizing geodesic from 
$\ov{\gamma}(\tau_0)$ to $\ov{\gamma}(s)$ in $U$. 
The speed of $\ov{\sigma}_s$ is $|\ov{\sigma}_s'|_h = L(\ov{\sigma}_s)/(s - \tau_0)$. Now lift this curve to $M$ via $\sigma_s: [\tau_0, s] \to M$ given by $\sigma_s(\tau) = \big(\tau, \ov{\sigma}_s(\tau)\big)$. To show that  
$\sigma_s$ is timelike (in fact, a timelike geodesic), it suffices to show $|\ov{\sigma}_s'|_h < 1$.  Since 
$\gamma$ is a timelike curve, we must have $|\ov{\gamma}'(\tau)|_h < 1$ for all $\tau \in [\tau_0, \tau_f]$. Integrating yields $L(\ov{\gamma}|_{[\tau_0,s]}) <  s - \tau_0$. Therefore
\begin{equation}\label{timelike lift}
|\overline{\sigma}_s'|_h = \frac{L(\overline{\sigma}_s)}{s - \tau_0} \leq \frac{L(\overline{\gamma}|_{[\tau_0,s]})}{s - \tau_0} < 1.
\end{equation}
Therefore $\sigma_s$ is a future directed timelike curve between $\gamma(\tau_0)$ and $\gamma(s)$. Now we define the homotopy $\ov{\Gamma}: [\tau_0,\tau_f] \times [\tau_0,\tau_f] \to \Sigma$ between $\ov{\gamma}$ and $\ov{\sigma}_{\tau_f}$ via  
\[
\ov{\Gamma}(s,\tau) = \big(\ov{\sigma}_s * \ov{\gamma}|_{[s,\tau_f]}\big)(\tau) = \left\{
\begin{array}{ll}
    \ov{\sigma}_s(\tau) & \text{ for } \tau_0 \leq \tau \leq s \\
      \ov{\gamma}|_{[s,\tau_f]} & \text{ for } s \leq \tau \leq \tau_f \\
   
\end{array} 
\right.
\]
and define $\Gamma:[\tau_0,\tau_f] \times [\tau_0, \tau_f] \to M$ by $\Gamma(s,\tau) = \big(\tau, \ov{\Gamma}(s,\tau)\big)$. We have shown that for each $s$, $\Gamma(s, \cdot)$ is a future directed timelike curve and $\Gamma(\tau_0, \cdot) = \gamma$ and $\Gamma(\tau_f, \cdot) = \sigma_{\tau_f}$. Thus $\Gamma$ is a future directed timelike homotopy between $\gamma$ and $\sigma_{\tau_f}$. 
\qed

\medskip

\begin{cor}\label{future connected cor}
Suppose at every point $p \in \Sigma$ there exists $0 \in U_p \subset T_p\S$
such that the exponential map, $\exp_p: U_p \to \S$, is a diffeomorphism onto $\Sigma$. Then any spacetime conformal to $(M, -d\tau^2 + h)$ is future one-connected. Hence $(M,g)$ is future one-connected.
\end{cor}

\medskip

\noindent{\bf Remark.}
We note that Corollary \ref{future connected cor} applies in particular to the case that $\S$ is a {\it Hadamard space} (i.e.\ a simply connected Riemannian manifold with nonpositive sectional curvature), for which we know that the exponential map is a global diffeomorphism about every point.

\medskip

For future divergence we have the following proposition and corollary. The ideas in the proofs also appear in Sbierski's paper \cite{Sbierski}.

\begin{prop}\label{future divergence prop}
Suppose $(\Sigma, h)$ is a complete Riemannian manifold and $I = (t_1, \infty)$, such that  $\t(t) \to \infty$ as $t \to \infty$.
Then $(M, g = -d\tau^2 + h)$ is future divergent. 
\end{prop}

\proof
Let $\gamma: [\tau_0, \infty) \rightarrow (M,g)$ be a future directed inextendible timelike curve parameterized by $\tau$: 
\[
\gamma(\tau) = \big(\tau, \overline{\gamma}(\tau)\big).
\]
Fix $T \in (\tau_0, \infty)$. Let $\overline{\sigma}: [\tau_0, T] \rightarrow \Sigma$ be a length minimizing geodesic between $\overline{\gamma}(\tau_0)$ and $\overline{\gamma}(T)$. Since $\overline{\sigma}$ is parameterized by $\tau$, the argument which led to (\ref{timelike lift}) also gives 
\[|\overline{\sigma}'|_h = L(\overline{\sigma})/(T - \tau_0) < 1.\]
Define $\sigma: [\tau_0, T] \rightarrow \Sigma$ by $\sigma(\tau) = \big(\tau, \overline{\sigma}(\tau)\big)$. Since $|\ov{\sigma}'|_h < 1$, $\sigma$ is a timelike curve (in fact timelike geodesic) connecting $\gamma(\tau_0)$ to $\gamma(T)$. Thus we have 
\begin{align}
d_g\big(\gamma(\tau_0), \gamma(T)\big) &\geq L_g(\sigma) \notag
\\
&= \int_{\tau_0}^T \sqrt{1 - |\overline{\sigma}'|_h^2}d\tau \notag
\\
&= \sqrt{(T - \tau_0)^2 - (T- \tau_0)^2|\overline{\sigma}'|_h^2} \notag
\\
&= \sqrt{(T - \tau_0)^2 - d_h^2\big(\overline{\gamma}(\tau_0), \overline{\gamma}(T)\big)} \label{distance}
\end{align}
where $d_h$ is the Riemannian distance function on $\Sigma$. 

Fix $\tau_1 \in [\tau_0, T]$. Then since $\gamma$ is timelike, we have $d_h\big(\overline{\gamma}(\tau_0), \overline{\gamma}(\tau_1)\big) < \tau_1 - \tau_0$ and $d_h\big(\overline{\gamma}(\tau_1), \overline{\gamma}(T)\big) < T - \tau_1$. Therefore there exists an $\epsilon > 0$ such that $d_h\big(\overline{\gamma}(\tau_0), \overline{\gamma}(\tau_1)\big) = \tau_1 - \tau_0 - \epsilon$. By the triangle inequality we have 
\begin{align*}
d_h\big(\overline{\gamma}(\tau_0), \overline{\gamma}(T)\big) &\leq d_h\big(\overline{\gamma}(\tau_0), \overline{\gamma}(\tau_1)\big) + d_h\big(\overline{\gamma}(\tau_1), \overline{\gamma}(T)\big)
\\
&< (\tau_1 - \tau_0 - \epsilon) + (T - \tau_1)
\\
&= T - \tau_0 - \epsilon.
\end{align*}
Using this in (\ref{distance}), we have
\begin{align*}
d\big(\gamma(\tau_0), \gamma(T)\big) &\geq \sqrt{(T - \tau_0)^2 - (T - \tau_0 - \epsilon)^2}
\\
&= \sqrt{2\epsilon(T - \tau_0) - \epsilon^2}.
\end{align*}
Therefore $\lim_{T \rightarrow \infty}d\big(\gamma(\tau_0), \gamma(T)\big) = \infty$.
\qed

\medskip
\begin{cor}\label{future divergence cor}
Assume the hypotheses of Proposition \ref{future divergence prop}. Then $(M,g)$ is future divergent so long as $a(t)$ is bounded away from $0$ for all large $t$.
\end{cor}

\proof
Let $\gamma: [\tau_0, \infty) \to M$ be a timelike curve in $(M,g)$ parameterized by $\tau$. Then 
\[
g\big(\gamma'(\tau),\gamma'(\tau)\big) = -|\gamma'(\tau)|_g^2 = a^2\big(t(\tau)\big)\big[-1 + |\ov{\gamma}'(\tau)|_h^2\big].
\]
Since $a(t)$ is bounded away from 0 for all large $t$, there exist $\tau_1 \in [\tau_0, \infty)$ and $b > 0$ such that $a\big(t(\tau)\big) > b$ for all $\tau \ge \tau_1$. So for all $\tau > \tau_1$, we have $|\gamma'(\tau)|_g^2 > b\big( 1 - |\ov{\gamma}'(\tau)|_h^2 \big)$, from which it follows that $L_g(\g|_{[\t_1, \t]}) >  b\, L_{g_0}(\g|_{[\t_1, \t]})$, where  $g_0 = -d\t^2 + h$. The result then follows from Proposition \ref{future divergence prop}.
\qed

\subsection{$C^0$-inextendibility of Anti-de Sitter space.}

In this subsection we illustrate how future one-connectedness and future divergence can be used to show that anti-de Sitter space is $C^0$-inextendible. There are various equivalent definitions of anti-de Sitter space. The most useful for us is the one where the metric is conformal to (part) of the Einstein static universe (cf.\ Carroll \cite{Carroll}). 

\medskip

\begin{Def}  
{\rm
The $d+1$-dimensional {\bf \emph{anti-de Sitter space (adS)}} is the spacetime $(\R \times S^d_+, g)$,  where $g$ is given by 
\[
g = \frac{1}{\cos^2\chi}\big[-dt^2 + d\chi^2 + \sin^2\chi d\Omega_{d-1}^2 \big],
\]
and $(\chi, \omega) \in [0,\frac{\pi}2) \times S^{d-1}$ are spherical coordinates on the (open) hemisphere 
$S^d_+$.
}

%
%
\end{Def}

\medskip

\begin{thm}
Anti-de Sitter space is $C^0$-inextendible.
\end{thm}

\proof 
We show adS space is future one-connected and future divergent so that Theorem \ref{future bdd empty} implies $\pd^+M = \emptyset$. Reversing the time orientation then shows that $\pd^-M= \emptyset$ and so Proposition \ref{bdd nonempty} implies that no $C^0$-extension of adS space exists. 

Since the round hemisphere satisfies the exponential map property, Corollary \ref{future connected cor} implies that anti-de Sitter space is future one-connected.   To show that adS space is future divergent, let $\gamma:[0,t_f) \to (\R^{d+1},g)$ be a future directed inextendible timelike curve parameterzed by $t$ (by a time translation we can assume $\gamma$ begins at $t =0$). For each $t \in [0, t_f)$, we have $\g(t) = (t, \chi(t), \omega(t))$, where $\omega$ represents coordinates on $S^{d-1}$.  There are essentially two cases to consider: There exists $t_k \nearrow t_f$ such that (i) $\lim_{k \to \infty} \chi(t_k) < \pi/2$ or (ii) $\lim_{k \to \infty} \chi(t_k) = \pi/2$.  In either case, for $k$ sufficiently large, there exists a $t$-line segment $\s_k$ from some $p_k \in \d I^+(\g(0))$ to $\g(t_k)$  Moreover, one has $\lim_{k \to \infty} L(\s_k) = \infty$. In case (i) this follows from the fact that we must have $t_f = \infty$.  In case (ii) this  follows from the fact that $\lim_{k \to \infty} \cos(\chi(t_k)) = 0$, so that the conformal factor becomes arbitrarily large, and that for $k$ sufficiently large,  $t_k - t(p_k)$ is uniformly positive.   
Since for each $k$ there exists a null geodesic from $\g(0)$ to $p_k$,
we conclude that adS space is future divergent.\qed  

\medskip

\noindent{\bf Remark.}
The ideas in this proof can also be used to prove that de Sitter space is $C^0$-inextendible, as Sbierski points out in his applications and open questions section~\cite{Sbierski}. The future one-connectedness of de Sitter space follows from Proposition~\ref{future connected prop} since the projection of timelike curves in de Sitter onto the spacelike spheres cannot contain any antipodal points.

\section{Open FLRW Spacetimes and Examples of $C^0$-Extensions}\label{FLRW}
In this section we define what is meant by an \emph{open} Friedmann-Lema{\^\i}tre-Robertson-Walker (FLRW) spacetime and apply the results of section 2 to show that these spacetimes are future one-connected and future divergent. Therefore if a 
$C^0$-extension exists for these spacetimes, then by Theorems \ref{future bdd empty} and \ref{reg result}, we have $\pd^+M = \emptyset$ and $\pd^-M$ must be an achronal locally Lipschitz hypersurface. Then we look at particular examples of open FLRW spacetimes which admit $C^0$-extensions but not $C^2$-extensions. It will be seen that the regularity of $\pd^-M$ cannot be improved. That is, we should not expect $\pd^-M$ to be a $C^1$-hypersurface. 

\subsection{Open FLRW Spacetimes}

\begin{Def}\label{flrw}
{\rm
An {\bf \emph{open FLRW spacetime}} is a spacetime $(M,g)$ where $M = (0, \infty) \times \R^d$ with coordinates $(t,r,\omega) \in (0,\infty)\times (0,\infty)\times S^{d-1}$ and the metric $g$ is 
\[
g = \left\{
\begin{array}{ll}
      -dt^2 + a^2(t)\big[dr^2 + r^2d\Omega_{d-1}^2\big] & \text{{\bf \emph{(Euclidean)}} }\\
      &\\
      -dt^2 + a^2(t)\big[dr^2 + \sinh^2(r)d\Omega_{d-1}^2\big] & \text{{\bf \emph{(Hyperbolic)}} }
\end{array} 
\right.
\]
The function $a:(0,\infty) \to (0,\infty)$ is called the {\bf \emph{scale factor}} and for us we demand that it satisfy the following four requirements
\begin{enumerate}

\item[(1)] $a$ is smooth.

\item[(2)] $a(0) := \lim_{t \to 0^+}a(t) = 0$.

\item[(3)] $a(t)$ has sublinear growth (i.e. there exist constants $m > 0$ and $b \geq 0$ such that $a(t) \leq mt +b$)

\item[(4)] $a'(t) > 0$ for all $t$.

\end{enumerate}
}
\end{Def}

\noindent{\bf Remark.}
The conditions on the scale factor are natural and are satisfied by the classical open FLRW spacetimes.  
With only minor modifications to our presentation, assumption (3) could be replaced by the assumption that 
$(M,g)$ is a future asymptotically simple and de Sitter spacetime, as defined, for example, in \cite{GBeem, AndGal}. This includes FLRW spacetimes with positive cosmological constant. In the context of the next theorem, such spacetimes are also future divergent.

\medskip

\begin{thm}\label{open FLRW thm}
Let $(M,g)$ be an open FLRW spacetime. If $(M_\ext, g_\ext)$ is a 
$C^0$-extension of $(M,g)$, then $\pd^+M = \emptyset$ and $\pd^-M$ is an achronal topological hypersurface. 
\end{thm}

\proof
We need to show $(M,g)$ is future one-connected and future divergent. Then the result follows from Theorems \ref{future bdd empty} and \ref{reg result}. Corollary \ref{future connected cor} implies $(M,g)$ is future one-connected. To show that $(M,g)$ is future divergent, define $\tau = \int_1^t\frac{1}{a(s)}ds$. Then condition (3) of the above definition implies that the range of $\tau$ is an infinite interval. Therefore $(M,g)$ is future divergent by Corollary \ref{future divergence cor}. 
\qed

\medskip

\noindent{\bf Remark.}
Since $t$ is a time function for an open FLRW spacetime $(M,g)$, one can reparameterize a future-directed timelike curve by $t$. Therefore if $(M_\ext,g_\ext)$ is a $C^0$-extension of $(M,g)$, then $\pd^-M$ is reached by 
timelike curves whose $t$-parameter approaches $0$.

\medskip

Now we present two classes of 
open FLRW spacetimes where $C^0$ extensions (but not necessarily smooth extensions) can be found. 

\subsection{$C^0$-Extendable 2-Dimensional Spacetimes}\label{2-dimensional subsection}

Let $(M,g)$ be a two-dimensional open FLRW spacetime.
We will find a $C^0$-extension $(M_\ext, g_\ext)$ of $(M,g)$. Since $\det(g) = -a^2(t)$, the metric is degenerate at $t = 0$ (in the case $a(0) = 0$), so we can't use the coordinates $(t,x)$ to extend the metric. Better coordinates to use are 
\begin{align*}
\tilde{t}(t,x) &= \int_0^ta(s)ds
\\
\tilde{x}(t,x) &= x - \int_1^t \frac{1}{a(s)}ds.
\end{align*}
A simple calculation shows that the metric in these coordinates is given by 
\begin{equation}\label{eq 4.1}
g = 2d\tilde{t}d\tilde{x} + a^2\big(t(\tilde{t})\big)d\tilde{x}^2.
\end{equation}
In these coordinates, we have $\det(g) = -1$ for all $(\tilde{t},\tilde{x})$, so no degeneracy in the metric occurs at $t = \tilde{t} = 0$. This allows us to extend $(M,g)$ into a larger spacetime. Extend $a$ onto $(-\infty,0)$ while keeping it continuous, so that now $\tilde{t}$ is defined on all of $\R$. Take $M_\ext = \R \times \Sigma$ with coordinates $(\tilde{t},\tilde{x})$ and metric $g_\ext$ defined by eq. (\ref{eq 4.1}). Then $(M,g)$ embeds isometrically as a proper subset of $(M_\ext,g_\ext)$. In fact $(M,g)$ is isometric to $(M_\ext, g_\ext)|_{\tilde{t} > 0}$. In this example $\pd^-M$ is given by the null hypersurface $\tilde{t} = 0$.

The scalar curvature $\mathcal{R}$ of $(M,g)$ is given by 
$\mathcal{R} = 2\frac{a''(t)}{a(t)}.$ This gives us a plethora of examples where we have a $C^0$-extension but not a $C^2$-extension. For example, by taking $a(t) = \sqrt{t}$ we have $\mathcal{R} = -\frac{1}{2t^2}$, so $\mathcal{R} \to -\infty$ as $t\to 0^+$, which is of course, an obstruction to a $C^2$-extension since any $C^2$-extension would have a continuous, finite-valued scalar curvature at all points of the extension. 

\subsection{Milne-Like Spacetimes}
The Milne universe is the spacetime $M = (0, \infty) \times \R^d$ with metric 
\[
g = -dt^2 + t^2\big[dr^2 + \sinh^2(r)d\Omega_{d-1}^2 \big].
\]
$(M,g)$ admits a smooth extension through $t = 0$. To see this define $T = t\cosh(r)$ and $R = t\sinh(r)$. Then one finds $(M,g)$ is isometric to the future cone $I^+(0) = \{(T,R,\omega) \mid  T > R \geq 0\}$ in Minkowski space $(\R^{d+1}, -dT^2 + dR^2 + R^2 d\Omega_{d-1}^2)$. Thus Minkowski space is a smooth extension of the Milne universe. In this case $\pd^-M$ is the boundary of $I^+(0)$, which 
shows that the regularity of $\pd^-M$ cannot be improved from $C^{0,1}$. 

We now define a class of spacetimes we call Milne-like. These spacetimes will admit extensions similar to the Milne universe, but we will find cases where the extension is $C^0$ but not $C^2$.

\begin{Def}
{\rm
A hyperbolic FLRW spacetime $(M,g)$ will be called {\bf \emph{Milne-like}} if $M = (0,\infty) \times \R^d$ and 
\[
g = -dt^2 + a^2(t)\big[dr^2 + \sinh^2(r)d\Omega_{d-1}^2\big],
\]
and the scale factor satisfies the following additional conditions:

\begin{itemize}

\item [(1)] $a'(0):= \lim_{t \to 0^+}a'(t) = 1$

\item [(2)] $\int_0^1\frac{1}{a(t)}dt = \infty$.

\end{itemize}

\noindent Put $b(t) = \exp\left(\int \frac{1}{a}\right)$ so that $b/b' = a$. Then we also require
\begin{itemize}

\item [(3)] $b'(0) := \lim_{t \to 0^+}b'(t) \in(0,\infty)$.
\end{itemize}

}
\end{Def}

\noindent{\bf Remark.}
Conditions (1) and (2) are necessary for condition (3). Without (2), we would have $\lim_{t \to 0^+}b(t) \neq 0$ so that $b'(0) = \lim_{t \to 0^+}\frac{b(t)}{a(t)} = \infty$. Given (3), we can apply L'H{\^o}pital's rule 
\[
b'(0) = \lim_{t \to 0^+}\frac{b(t)}{a(t)}= \lim_{t \to 0^+}\frac{b'(t)}{a'(t)}, 
\]
from which we see that we must have $a'(0) = 1$.

\medskip

\medskip

\begin{thm}\label{Milne-like thm}
Let $(M,g)$ be any Milne-like spacetime. Then $(M,g)$ is $C^0$-extendible.
\end{thm}

\proof
Define $T = b(t)\cosh(r)$ and $R = b(t)\sinh(r)$. Then 
\begin{align}\label{Milne metric}
g &=-dt^2 + a^2(t)\big[dr^2 + \sinh^2(r)d\Omega_{d-1}^2\big] \notag
\\
&= \frac{1}{\big[b'\big(t(T,R)\big) \big]^2} \big[-dT^2 + dR^2 + R^2d\Omega_{d-1}^2\big].
\end{align}

Note that $b \to \infty$ as $t \to \infty$ since we are assuming $a$ has sublinear growth (condition (3) in the definition of open FLRW spacetime). Therefore $(M,g)$ is isometric to the region 
$\{(T,R,\omega) \mid  R < T < \infty\}$. 
Since $b'(0) \in (0,\infty)$, equation \ref{Milne metric} implies that there is no degeneracy in the metric as $t \to 0^+$ in the $(T,R,\omega)$ coordinate system, so we can extend the metric through $t = 0$. There are of course an infinite number of ways to do this. For specific choices of $a(t)$ certain extensions are more readily apparent. 
See the examples below. For our general scenario, we can extend by keeping $b'(0)$ constant through $t = 0$.
For this choice, our extended manifold is 
\[
M_\ext = \R^{d+1} =  \{(T,R,\omega) \mid R < T < \infty\} \cup \{(T,R,\omega) \mid -\infty < T \leq R\}
\]
and the extended metric is
\[   
g_\ext =\left\{
\begin{array}{ll}
      \frac{1}{\big[b'\big(t(T,R)\big) \big]^2} \big[-dT^2 + dR^2 + R^2d\Omega_{d-1}^2\big] & \text{ on } \{(T,R,\omega) \mid R < T < \infty\}\\
      \\
      \frac{1}{\big[b'(0) \big]^2} \big[-dT^2 + dR^2 + R^2d\Omega_{d-1}^2\big] & \text{ on } \{(T,R,\omega) \mid -\infty < T \leq R\} \\
\end{array} 
\right. \]
Then $(M_\ext, g_\ext)$ is a $C^0$ extension of $(M,g)$.
\qed

\medskip

Let's look at a couple of examples. 

\medskip

\begin{ex}
{\rm 
Consider $a(t) = \tanh(t)$. In this case we have $b(t) = \sinh(t)$ and so 
\begin{align}
g &= -dt^2 + \tanh^2(t)\big[dr^2 + \sinh^2r d\Omega_{d-1}^2\big] \notag
\\
 &= \frac{1}{\cosh^2(t)}\big[ -dT^2 + dR^2 + R^2d\Omega_{d-1}^2\big] \notag
 \\
 &= \frac{1}{1 + (T^2 - R^2)}\big[-dT^2 + dR^2 + R^2 d\Omega_{d-1}^2\big] 
\end{align}
By defining 
\begin{equation}
g_\ext = \frac{1}{1 + (T^2 - R^2)}\big[-dT^2 + dR^2 + R^2 d\Omega_{d-1}^2\big],
\end{equation}
and $M_\ext = \{(T,R,\omega) \mid 1 + T^2 - R^2 >0\}$, we find $(M_\ext, g_\ext)$ is a $C^\infty$-extension of $(M,g)$. 
}
\end{ex}

\medskip

\begin{ex}
{\rm
Now let's consider an example where we have a $C^0$-extension but not a $C^2$-extension. Let $a(t)$ be a function satisfying the conditions of Definition \ref{flrw} such that $a(t) = t + t^2$ for $t \leq 100$ and extend it past $t \geq 100$, so that it has sublinear growth and satisfies $a'>0$ (this is so that we satisfy the conditions of open FLRW spacetimes). For $t < 100$, let $b(t) = t/(1 + t)$ so $b/b' = a$. Since $b = \sqrt{T^2 - R^2}$, we have $t = \frac{\sqrt{T^2 - R^2}}{1 - \sqrt{T^2 - R^2}}$
We find
\begin{align}
g&= -dt^2 + (t+t^2)^2\big[dr^2 + \sinh^2rd\Omega_{d-1}^2 \big] \notag
\\
&= \frac{1}{\big[b'(t)\big]^2}\big[-dT^2 + dR^2 + R^2d\Omega_{d-1}^2 \big] \notag
\\
&= (1 + t)^4\big[-dT^2 + dR^2 + R^2d\Omega_{d-1}^2 \big] \notag
\\
&= \left(\frac{1}{1 - \sqrt{T^2 - R^2}} \right)^4\big[-dT^2 + dR^2 + R^2d\Omega_{d-1}^2 \big] \,.
\end{align}
To extend $(M,g)$, we can take 
\begin{equation}
g_\ext = \left(\frac{1}{1 - \sqrt{|T^2 - R^2|}} \right)^4\big[-dT^2 + dR^2 + R^2d\Omega_{d-1}^2 \big]
\end{equation}
and
\[
M_\ext = \{(T,R,\omega) \mid |T^2 - R^2| < 1\} \,.
\]
Then $(M_\ext, g_\ext)$ is a $C^0$-extension of $(M,g)$. There can be no $C^2$-extension of $(M,g)$ since the scalar curvature of $(M,g)$ is (see  \cite[p.\ 116]{BEE})  
\begin{equation}
\mathcal{R} =\frac{-d(d-1)}{a^2(t)} + 2d\left(\frac{a'(t)}{a(t)}\right)^2 + 2d\frac{a''(t)}{a(t)} + (d^2 - 3d)\left(\frac{a'(t)}{a(t)}\right)^2
\end{equation} 
and so we find
\[
\lim_{t \to 0^+}\mathcal{R}= \lim_{t \to 0^+}2d\frac{a''(t)}{a(t)} = \lim_{t \to 0^+}\frac{4d}{t + t^2} = \infty.
\]
}
\end{ex}

\medskip
\noindent{\bf Remark.} This example satisfies the weak energy condition for small $t$. There are also Milne-like examples that satisfy the strong energy condition.  Howevever, it can be seen that no Milne-like spacetime can satisfy both the weak and strong energy conditions unless it is exactly the Milne spacetime.

\section{Spherical Symmetry}\label{sphsym}

The open FLRW spacetimes possess spherical symmetry, so we wish to describe this spherical symmetry in the class of $C^0$ spacetimes. The definition of spherical symmetry given in \cite[Box 23.3]{MTW}, makes sense at the $C^0$ level.  There, it is assumed that the group of isometries of spacetime $(M^{d+1},g)$ contains $SO(d)$ as a subgroup, such that the orbits of this action are spacelike $(d-1)$-spheres ($d = 3$ in their discussion).  It is further assume that there exists a timelike vector field $u$ 
 invariant under the $SO(d)$ group action. Then, under these assumptions, their arguments lead (in the case $d =3$; in fact any $d$ odd would suffice) to the existence about every point of $M$ local coordinates $(x,y,\omega \in S^{d-1})$ such that 
with respect to these coordinates
the metric takes the form 
\begin{equation}\label{spherically symmetric metric}
g = A(x,y)dx^2 + 2B(x,y)dxdy + C(x,y)dy^2 + R^2(x,y)d\Omega_{d-1}^2  \,.
\end{equation}
If coordinates $(x,y, \omega)$ can be introduced so that the metric takes this form, we will say that spacetime is {\bf \emph{spherically symmetric}} and will refer to the coordinates 
$(x,y,\omega)$ as {\bf \emph{spherically symmetric coordinates}}.  The choice of radial function $R$ is unique in the following sense: If $(x,y,\omega)$ and $(\bar{x}, \bar{y},\omega)$ are spherically symmetric coordinates, such that $x$ and $y$ are solely functions of $\bar{x}$ and $\bar{y}$, then both coordinate systems induce the same radial function on the overlap.  It should be noted that the usual procedure one uses to eliminate the cross term cannot be applied in the $C^0$ setting because this requires a Lipschitz condition on $A$, $B$, and $C$.

We will say that $(M,g)$ is {\bf \emph{strongly spherically symmetric}} if about every point there are coordinates $(T,R,\omega)$ such that in this coordinate neighborhood the metric takes the form 
\begin{equation}\label{strong ssm}
g = -F(T,R)dT^2 + G(T,R)dR^2 + R^2d\Omega_{d-1}^2  \,,
\end{equation}
and we call $(T,R,\omega)$ {\bf \emph{strongly spherically symmetric coordinates}}.  To achieve the metric form \eqref{strong ssm} via a change of coordinates from \eqref{spherically symmetric metric}, requires greater regularity on the metric, at least $C^1$, and, in addition, a $C^1$ genericity condition on $R$.
We note that Milne-like spacetimes are strongly spherically symmetric if one defines the radial function to be $\wt{R} = R/b'$.  

We were able to find $C^0$-extensions for a Milne-like spacetime $(M,g)$ by writing $(M,g)$ in strongly spherically symmetric coordinates. A natural question to ask is: Can strongly spherically symmetric coordinates be used to find a $C^0$-extension for Euclidean FLRW spacetimes? What about hyperbolic FLRW spacetimes that are not Milne-like? The results in the following subsections answer in the negative.

\medskip
\subsection{No Strongly Spherically Symmetric Extensions for Euclidean FLRW Spacetimes}

\begin{thm}\label{euclid no S ext}  Let $(M,g)$ be a Euclidean FLRW spacetime where the scale factor $a(t)$ satisfies $a'(0) := \lim_{t \to 0^+}a'(t) \in (0,\infty]$.
Then, subject to a suitable initial condition, there exists a unique transformation of the form,
\begin{equation}\label{T and R coords}
T = T(t,r)  \quad R = R(t,r)  
\end{equation}
such that $g$ takes the strongly spherically symmetric form 
$$\label{sss}
g = -F(T,R)dT^2 + G(T,R)dR^2 + R^2d\Omega_{d-1}^2,
$$
where $F$ and $G$ are regular (away from a curve in the $r$-$t$ plane along which the Jacobian determinant $J(r,t) = \frac{\d (T,R)}{\d (t,r)}$ vanishes). 

Now suppose that $M$ admits a $C^0$-extension $M_{\ext}$, and consider the behavior of the metric in these coordinates
on approach  to $\pd^-M$ (cf.\ Theorem \ref{open FLRW thm}).  Let $\g:[0,1] \to M_{\ext}$ be a future directed timelike curve with past end point $\g(0) \in \pd^-M$, and suppose $R$ has a finite positive limit along $\g$ as $t \to 0^+$. (Note, by the achronality of $\pd^-M$, 
$\g((0,1]) \subset M$.)   Then the following hold along~$\g$.

\begin{enumerate}

\item[(a)] $\lim_{t \to 0^+} G = 0$.

\item[(b)] If $F$ has a finite nonzero limit as $t \to 0^+$, then $T \to \pm \infty$ as $t \to 0^+$. 

\end{enumerate}

\end{thm}

\noindent\emph{Remark.} By a `suitable initial condition', we mean the following: The transformation \ref{T and R coords} is unique up to a function $f$ which is determined by specifying $T$ along a certain curve in the first quadrant of the $(r,t)$-plane. This is shown in the proof below.

\proof We begin by solving explicitly for $R$, $T$, $G$, and $F$ in terms of 
$t$, and $r$.
Immediately, we find $R = ra(t)$. To see this, consider a codimension 2 surface of constant $T$ and $R$. Since $T$ and $R$ are functions of $t$ and $r$ only, a surface of constant $T$ and $R$ corresponds to a surface of constant $t$ and $r$. By restricting the metric to this surface, we have $R^2d\Omega_{d-1}^2  = r^2a^2(t)d\Omega_{d-1}^2$ and hence $R = ra(t)$. 

Let $T_t = \pd T/ \pd t$ and $T_r = \pd T/ \pd r$. Then 

\begin{align*}
dT^2 &= T_t^2dt^2 + 2T_tT_rdtdr + T_r^2dr^2
\\
dR^2 &= r^2a'^2 dt^2 +2raa'dtdr + a^2dr^2.
\end{align*}
So we want
\begin{equation}\label{metric match}
-dt^2 + a^2(t)\big[dr^2 + r^2d\Omega_{d-1}^2\big] = -FdT^2 + GdR^2 + R^2d\Omega_{d-1}^2
\end{equation}
From equation (\ref{metric match}), we find

\begin{align} 
-1 &= -FT_t^2 + Gr^2 a'^2 &&\Longrightarrow  &&&FT_t^2= 1 + Gr^2a'^2 \label{finding F and G 1}
\\
0 &= -FT_tT_r + Graa' &&\Longrightarrow &&&F^2T_t^2T_r^2 = G^2r^2a^2a'^2 \label{finding F and G 2}
\\
a^2 &= -FT_r^2 + Ga^2  &&\Longrightarrow &&&FT_r^2 = a^2(G-1) \label{finding F and G 3}
\end{align}
By substituting (\ref{finding F and G 1}) and (\ref{finding F and G 3}) into (\ref{finding F and G 2}),  we find 

\begin{equation}\label{G}
G(r,t) = \frac{1}{1 - r^2a'^2}
\end{equation}
Substituting this into (\ref{finding F and G 1}) and (\ref{finding F and G 3}), we find 
\begin{equation}\label{eq for T_t}
FT_t^2 = \frac{1}{1 - r^2a'^2} \:\:\:\: \text{ and } \:\:\:\: FT_r^2 = \frac{r^2a^2a'^2}{1-r^2a'^2}
\end{equation}
Therefore $(T_r/T_t)^2 = (raa')^2$. Since we require the metric to be Lorentzian, the leftmost equation in (\ref{finding F and G 2}) implies that we must have $T_r/T_t = raa'$. A solution to this PDE must be constant along the integral curves of $dt/dr = -raa'$ in the $(t,r)$-plane, so a general solution for $T$ is 
\begin{equation}\label{Teq}
T(r,t) = f\left(\frac{r^2}{2} + \int\frac{1}{aa'}\right)
\end{equation}
where $f$ is some smooth function. $f$ is uniquely determined by specifying $T$ on a curve which is transversal to the curves $\frac{r^2}{2} + \int\frac{1}{aa'} = \text{const.}$ Thus there is a degree of freedom when choosing strongly spherically symmetric coordinates.

In summary we have 
\begin{itemize}
\item $R = ra(t)$ 

\item $T = f\left(\frac{r^2}{2} + \int \frac{1}{aa'}\right)$

\item $G = \frac{1}{1 - r^2a'^2}$

\item $F = GT_t^{-2} = G\left(\frac{aa'}{f'}\right)^2.$
\end{itemize}

\noindent The Jacobian of the transformation is 
\[
J = T_rR_t - T_tR_r = f'[r^2 a' - 1/a']. 
\]
Therefore $F$ and $G$ are regular everywhere except where the Jacobian vanishes, namely along the curve $r^2a'(t)^2 = 1$ (since, from \eqref{eq for T_t} and \eqref{Teq}, $f' \ne 0$).  Also, note that $T$ and $R$ change causal character here. 

We can write the metric as
\begin{align}
g &= -FdT^2 + GdR^2 + R^2d\Omega^2 \notag
\\
 &= \frac{1}{1 -r^2a'^2}\left[-\left(\frac{aa'}{f'}\right)^2dT^2 + dR^2  \right] + R^2 d\Omega^2 \label{g in r,t,T,R} \,.
\end{align}
In equation (\ref{g in r,t,T,R}), $r$ and $t$ are smooth implicit functions of $R$ and $T$ away from $r^2a'(t)^2 = 1$.

Now restrict to $\g$. Along $\g$
we see $G \to 0$ as $t \to 0^+$ since $r = R/a$. This establishes (a).
To prove (b), let us use $s$ to denote the argument of $f$. Then
\begin{align}
s(R,t) &= \frac{R^2}{2a^2} + \int \frac{1}{aa'}\notag
\\
&= \frac{\frac{1}{2}R^2 + a^2\int \frac{1}{aa'}}{a^2} \label{s}
\end{align}
$F$ is given by  
\begin{align}
F(t) &= \frac{1}{1 - R^2(a'/a)^2}\left(\frac{aa'}{f'\big(s(t)\big)} \right)^2 \notag
\\
&= \left(\frac{a^4}{(a/a')^2 - R^2}\right)\frac{1}{\big[f'\big(s(t)\big)]^2} \label{F}
\end{align}
Rearranging (\ref{s}) and (\ref{F}) gives us
\begin{align}
s^2\big[f'(s)\big]^2 &= \left(\frac{\frac{1}{2}R^2 + a^2\int \frac{1}{aa'}}{a^2}\right)^2\left(\frac{a^4}{F\big[(a/a')^2 - R^2\big]} \right) \notag
\\
&= \frac{\big[ \frac{1}{2}R^2 + a^2\int \frac{1}{aa'}\big]^2}{F\big[(a/a')^2 - R^2\big]} \label{sf'}
\end{align}
Now assume $F$ has a finite  nonzero limit as $t \to 0^+$.   
Since $(a/a') \to 0$ and $a^2\int\frac{1}{aa'} \to 0$ (by L'H{\^o}pital's rule) as $t \to 0^+$ along $\gamma$,
there is a constant 
$0 < c < \infty$ such that
\[
\lim_{t \to 0^+}s^2\big[f'(s)\big]^2 = c^2
\]
Note that $t \to 0^+$ implies $s \to \infty$ along $\gamma$. Therefore 
the above limit is equivalent to $\lim_{s \to \infty}s^2\big[f'(s)\big]^2 = c^2$. 
As noted above, $f' \ne 0$,
so it follows that $\lim_{s \to \infty}sf'(s) = \pm c$. Fix $0 < \epsilon < c/2$. Then there exists an $S$ such that $s > S$ implies $|sf'(s) \mp c| < c/2$ so $f'(s) > \pm\frac{c}{2s}$. By integrating over all $s > S$, we find that $f(s) \to \pm \infty$ as $s \to \infty$. Hence $T \to \pm \infty$ as $t \to 0^+$ along $\gamma$.
\qed

\medskip

\begin{cor}\label{euclid no R ext}
Let $(M,g)$ be a Euclidean FLRW spacetime where the scale factor satisfies 
the condition in Theorem \ref{euclid no S ext}. Then there is no $C^0$ strongly spherically symmetric extension of $(M,g)$.   
\end{cor}

\noindent {Remark.} By the conclusion we mean precisely the following: There is no point $p \in \pd^-M$ for which there exist strongly spherically symmetric coordinates $(T,R, \omega)$ defined on a neighborhood $U$ of 
$p$ such that on  $U \cap M$, $T$ and $R$ are functions of $t$ and $r$ only and $g = \psi^*g_{ext}$, where $\psi$ is the transformation: $(t,r) \to (T,R)$.

\proof  Suppose there is such an extension. Then there is a point $p \in \pd^-M$ and a neighborhood $U$ of $p$ with strongly spherically symmetric coordinates $(T, R, \omega)$ such that $T$ and $R$ are as in the remark.  
But then the conclusions (a) and (b)  of Theorem~\ref{euclid no S ext} apply. In particular, Theorem \ref{euclid no S ext} implies that $G(p) = 0$, so the extended metric is degenerate at $p$.
\qed

 
%
%
%

\subsection{No Strongly Spherically Symmetric Extensions for non-Milne-like Hyperbolic FLRW Spacetimes}

We have analogous statements of Theorem \ref{euclid no S ext} and Corollary \ref{euclid no R ext} 
for hyperbolic FLRW spacetimes. However, we have to rule out the Milne-like spacetimes since we know these admit $C^0$ extensions by Theorem \ref{Milne-like thm}.  

\medskip

\begin{thm}\label{hyper no S ext}  Let $(M,g)$ be a  hyperbolic FLRW spacetime where the scale factor $a(t)$ satisfies 
\[
a'(0) := \lim_{t \to 0^+} a'(t) \in[0,\infty], \quad a'(0) \neq 1  \,.
\]
Then, subject to a suitable initial condition, there exists a unique transformation of the form,
\[
T = T(t,r)  \quad R = R(t,r)  
\]
such that $g$ takes the strongly spherically symmetric form 
\[
g = -F(T,R)dT^2 + G(T,R)dR^2 + R^2d\Omega_{d-1}^2,
\]
where $F$ and $G$ are regular (away from a curve in the $r$-$t$ plane along which the Jacobian determinant $J(r,t) = \frac{\d (T,R)}{\d (t,r)}$ vanishes).   Suppose $M$ admits a $C^0$-extension 
$M_{\ext}$.
Let $\g$ be a timelike curve in $M_{\ext}$ with past end point on $\pd^-M$, such that $R$ has a finite positive limit  along $\g$ as $t \to 0^+$.
Then we have $\lim_{t \to 0^+}G = 0$ along $\g$.

\end{thm}

\medskip 

\proof
The proof is hardly different from the proof of Theorem \ref{euclid no S ext}. The same analysis leads to the following expressions for $R$, $T$, $G$, and $F$
 \begin{itemize}
\item $R = \sinh(r)a(t)$ 

\item $T = f\bigg(\ln\big(\cosh(r)\big) + \int \frac{1}{aa'}\bigg)$

\item $G = \big[\cosh^2(r) - \sinh^2(r)a'^2(t)\big]^{-1}$

\item $F = \cosh^2(r)GT_t^{-2} = \cosh^2(r)G\left(\frac{aa'}{f'}\right)^2.$
\end{itemize}
where $f$ is some differentiable function which is uniquely determined by specifying $T$ on a curve which is transversal to the curves $\ln\big(\cosh(r)\big) + \int\frac{1}{aa'} = \text{const.}$

We have
\begin{align}
G &= \big[\cosh^2\big(\sinh^{-1}(R/a) \big) - \sinh^2\big(\sinh^{-1}(R/a) \big)a'^2\big]^{-1} \notag
\\
&=\big[(R/a)^2 + 1 - R^2a'^2/a^2\big]^{-1}\notag
\\
&= \frac{a^2}{R^2(1 - a'^2) + a^2}
\end{align}
Since $a'(0) \neq 1$, it follows that $G \to 0$ as $t \to 0^+$ along $\g$.
\qed 

\medskip


\medskip

\begin{cor}\label{hyper no R ext}
Let $(M,g)$ be a hyperbolic FLRW spacetime where the scale factor satisfies the condition in Theorem \ref{hyper no S ext}. Then there is no $C^0$ strongly spherically symmetric extension of $(M,g)$.   

\proof The remark following Corollary \ref{euclid no R ext} still applies.  The proof is then essentially the same as the proof of  Corollary \ref{euclid no R ext}.

\end{cor}


\medskip

The theorems and corollaries in this and the previous  subsection show that if $C^0$ extensions exist for open FLRW spacetimes (which are not Milne-like), then the extensions are likely not to be spherically symmetric. However, one might speculate 
that if an extension exists, then a spherically symmetric extension should exist as well. Thus we propose, with rather limited evidence, the following two conjectures. If these two conjectures are true, they would prove that open FLRW spacetimes (which are not Milne-like), are $C^0$-inextendible. 

\medskip
\noindent {\it Conjecture 1.} If an open FLRW spacetime admits a $C^0$ extension, then it admits a spherically symmetric $C^0$ extension. 

\medskip

\noindent {\it Conjecture 2.} An open FLRW spacetime (which is not Milne-like) admits no spherically symmetric $C^0$ extensions.

\medskip

For a quite different analysis of FLRW models near the big bang, we mention the
paper \cite{Klein}.  The $3+1$ dimensional models considered there exhibit degeneracy of the metric at the big bang in the Fermi-coordinates of a co-moving observer considered by the authors, and thus do not lead to $C^0$ extensions in the sense of the present paper.  

\section{Further Remarks on $\pd^-M$}
In this section we make some brief remarks on the structure of $\pd^-M$. For $d > 1$, we will see that the structure of $\pd^-M$ limits the possible extensions one can find. We will also show that these limitations do not exist when $d = 1$ (i.e. two-dimensional spacetimes). 

Consider a spacetime  $M = (0,\infty) \times \Sigma$, $g = -dt^2 + a^2(t)h$, with $(\Sigma, h)$  a Riemannian manifold of dimension $d >1$, such that $(M,g)$ is future divergent and future one-connected (an FLRW spacetime, for example).  Suppose $(M_\ext, g_\ext)$ is a $C^0$ extension of $(M,g)$.  By Theorem~\ref{reg result},
$\pd^-M $ is an achronal $C^0$ hypersurface.   It can be represented locally as a graph over a smooth hypersurface, where the graphing function satisfies a Lipschitz condition.  Hence, as a consequence of Rademacher's theorem,  $\pd^-M$ is differentiable almost everywhere, in the sense of having a well-defined tangent plane at almost all points.   Using the continuity of $g_\ext$ and the achronality of $\pd^-M$, one can show that  these tangent planes cannot be timelike; tangent vectors to $\pd^-M$, when they exist are spacelike or null.

Consider a point $p \in \pd^-M$ and let $\{y^0, y^1, \dotsc, y^d\}$ be coordinates for a neighborhood $U$ of $p$ with $\pd/ \pd y^0$ timelike and future directed. We can find a line $\sigma: [0,1] \to \R^d$, 
\[
\sigma(s) = \big(a^1s, \dotsc, a^ds\big), \text{ where } a^1, \dotsc, a^d \text{ are constants},
\]
such that the timelike surface
\[
T = \big\{\big(y^0, \sigma(s)\big) \mid s \in [0,1]\big\} \cap U
\]
intersects $\pd^-M$ in a curve $c:[0,1] \to \pd^-M$, $s \to c(s)$,  which is necessarily achronal in $T$.  Hence $c$ is differentiable a.e., with tangent vectors that are spacelike or null.  By suitably adjusting $T$ one can ensure that $c$ is differentiable and spacelike at some point. Let's assume in fact that we can choose $c$ to be spacelike on a set of positive of measure.  

We define $c_n: [0,1] \to T$ by intersecting $T$ with $\{t = 1/n\}$. We have $c_n(s) = \big(n^{-1}, \ov{c}_n(s)\big) \in (0,\infty) \times \Sigma$, where $\ov{c}_n$ is a spacelike curve in $(\Sigma, h)$. By Dini's theorem $c_n$ converges uniformly to $c$.   Let us make the assumption that  $c_n'$ converges to $c'$ a.e. This would be true (for a subsequence) by Arzel{\'a}-Ascoli if $\{c_n'\}$ satisfied a H{\"o}lder condition.  This, of course, cannot be expected to hold in general, but does hold in our two-dimensional examples.  
If we make this assumption, then basic analysis implies that the length of $c_n$ converges to that of $c$, which itself has positive length,
\begin{equation} \label{length of c eq}
\lim_{n \to \infty}\int_0^1\sqrt{g_\ext \big( c_n'(s),c_n'(s)\big)}ds = \int_0^1 \sqrt{g_\ext\big( c'(s),c'(s)\big)}ds > 0 \,.
\end{equation}
Now note that  
\begin{equation}\label{length eq}
\int_0^1\sqrt{g_\ext \big( c_n'(s),c_n'(s)\big)}ds = a(n^{-1})\int_0^1 \sqrt{h\big(\ov{c}_n'(s),\ov{c}_n'(s)\big)}ds \,.
\end{equation}
Let $(x^1, \dotsc, x^d)$ be coordinates for $\Sigma$ and set $x^i_n = x^i \circ c_n$. Then we have 
\begin{equation}\label{h eq}
h\big(\ov{c}_n'(s),\ov{c}_n'(s)\big) = h_{ij}\frac{dx^i_n}{ds}\frac{dx^j_n}{ds} \:\:\:\:\:\: i,j = 1, \dotsc, d \,.
\end{equation}
However, we can induce coordinates $(y^1_n, \dotsc, y^d_n)$ on $\{n^{-1}\} \times \Sigma$ by the graphing function $y^0 = y^0(y^1, \dotsc, y^d)$. Then the chain rule gives 
\begin{equation}
\frac{dx^i_n}{ds} = \frac{\pd x^i}{\pd y^k_n}\frac{d y^k_n}{ds} = a^k\frac{\pd x^i}{\pd y^k} \:\:\:\:\:\: i,k = 1, \dotsc, d \,.
\end{equation} 
Therefore equation (\ref{h eq}) gives 
\begin{equation}
h\big(\ov{c}_n'(s),\ov{c}_n'(s)\big) = h_{ij} a^ka^l\frac{\pd x^i}{\pd y^k_n}\frac{\pd x^j}{\pd y^l_n}.
\end{equation}
The right hand side of (\ref{length eq}) then becomes,
\begin{equation} \label{length in terms of y eq}
a(n^{-1})a^ka^l\int_0^1 \sqrt{ h_{ij}\frac{\pd x^i}{\pd y^k_n}\frac{\pd x^j}{\pd y^l_n}}ds \,.
\end{equation}

Let's suppose we are dealing with a Euclidean FLRW model so that $h_{ij} = \delta_{ij}$. Then  
a contradiction to equation (\ref{length of c eq}) will result if  the following condition holds:
\begin{equation}\label{blowup}
\lim_{n \to \infty} a(n^{-1})\sup_{s \in [0,1], i,k}\left\{\left|\frac{\pd x^i}{\pd y_n^k}\big(c_n(s)\big) \right| \right\} = 0.
\end{equation}
Hence, 
this is telling us  that, in order for an extension to exist, the coordinates on $\Sigma$ and the coordinates about $p \in \pd^-M$ must behave in a certain asymptotic manner as one approaches $\pd^-M$: Roughly speaking, some of the partial derivatives appearing in \eqref{blowup} must become unboundedly large on approach to $\pd^-M$, at a rate based on the rate at which $a(t) \to 0$. 

We don't see this asymptotic behavior in our two-dimensional spacetimes which admitted $C^0$ extensions that we explored in Subsection \ref{2-dimensional subsection}. Indeed in those spacetimes we have $\pd x /\pd \tilde{x} = 1$.   A contradiction is avoided in this case because $\pd^-M$ is a {\it one}-dimensional null hypersurface, and hence does not admit any spacelike directions.  

\newpage

\bibliographystyle{amsplain}
\bibliography{extension}

\end{document}